# SAMPLED-DATA BOUNDARY FEEDBACK CONTROL OF 1-D HYPERBOLIC PDES WITH NON-LOCAL TERMS

**Iasson Karafyllis[*] and Miroslav Krstic[**]**

[*]Dept. of Mathematics, National Technical University of Athens, Zografou Campus, 15780, Athens, Greece, email: iasonkar@central.ntua.gr

[**]Dept. of Mechanical and Aerospace Eng., University of California, San Diego, La Jolla, CA 92093-0411, U.S.A., email: krstic@ucsd.edu

**Abstract**

The paper provides results for the application of boundary feedback control with Zero-Order-Hold (ZOH) to 1-D linear, first-order, hyperbolic systems with non-local terms on bounded domains. It is shown that the emulation design based on the recently proposed continuous-time, boundary feedback, designed by means of backstepping, guarantees closed-loop exponential stability, provided that the sampling period is sufficiently small. It is also shown that, contrary to the parabolic case, a smaller sampling period implies a faster convergence rate with no upper bound for the achieved convergence rate. The obtained results provide stability estimates for the sup-norm of the state and robustness with respect to perturbations of the sampling schedule is guaranteed.

## 1. Introduction

Sampled-data feedback control has been studied extensively for finite-dimensional systems due to the use of digital technology in modern control systems for the implementation of the controller (see for instance [7,11,13,15,23,24,33] and the references therein). However, sampled-data feedback control has been scarcely studied for infinite-dimensional systems. Most of the available results deal with delay systems (see [8,14,26,27,28,29,32]). For systems described by Partial Differential Equations (PDEs) the design of sampled-data feedback control faces major technical difficulties: even the notion of the solution of a PDE under sampled-data feedback control has to be clarified. Sampled-data controllers for parabolic systems were designed by Fridman and coworkers in [1,9,10,31] by using matrix inequalities. The works [21,30] provided necessary and sufficient conditions for sampled-data control of general infinite-dimensional systems under periodic sampling (see also [22,35] for the case of "generalized sampling"). Approximate models of infinite-dimensional systems were used in [34] for practical stabilization. A sampled-data feedback controller for hyperbolic age-structured models was proposed in [17].

In the linear finite-dimensional case, there are results that guarantee closed-loop exponential stability for continuous-time linear feedback designs when applied with Zero-Order-Hold (ZOH) and sufficiently small sampling period (see for instance [11,13,23,24]). The results deal with the globally Lipschitz case (which contains the linear case as a subcase) and the application of the continuous-time feedback under ZOH is called the "emulation" sampled-data feedback design.

A general robustness result that guarantees closed-loop exponential stability for continuous-time linear boundary feedback designs when applied with ZOH and arbitrary (not necessarily periodic) sampling schedules of sufficiently small sampling period is missing for the case of systems described by PDEs. In the recent work [18], efforts were made towards the development of such general results for linear parabolic PDE systems.

While the development of continuous-time boundary feedback controllers for hyperbolic PDE systems has progressed significantly during the last decade (see [2,3,5,16,19,20] for a single PDE and [4,6,12] for systems of PDEs), there are no results that guarantee stability properties for the sample-and-hold implementation of continuous-time controllers with arbitrary sampling schedules of sufficiently small sampling period. The present paper provides sampled-data, boundary feedback controllers for 1-D, first-order, linear, hyperbolic PDEs with non-local terms. The design is based on the emulation of the continuous-time boundary feedback design presented in [19]. It is proved that there is a sufficiently small sampling period, such that the closed-loop system preserves exponential stability under the sample-and-hold implementation of the controller (Theorem 2.1). In other words, we prove that emulation design works for the case of linear hyperbolic PDEs with boundary feedback. The derived exponential stability estimates are expressed in the sup-norm of the state and (conservative) upper bounds for the sampling period are derived. Finally, robustness with respect to the sampling schedule is established, exactly as in the finite-dimensional case.

The methodology for proving the main result of the present work is very different from the corresponding methodology in the parabolic case. While both proofs of the main results in [18] exploit an eigenfunction expansion procedure, the proof of Theorem 2.1 relies on the representation of the solution of the closed-loop system by means of an Integral Delay Equation (IDE), as proposed in [16]. However, there is an additional important difference between the parabolic and the hyperbolic case. In the hyperbolic case (Theorem 2.1), by selecting a sufficiently small maximum allowable sampling period we can achieve an arbitrarily fast rate of convergence. This is not possible in the parabolic case. This important difference can be explained by the fact that the nominal continuous-time feedback law (proposed in [19]) achieves finite-time stability in the hyperbolic case, while the nominal continuous-time feedback laws in the parabolic case achieve exponential stability. The proof of Theorem 2.1 provides an estimate of how small the maximum allowable sampling period must be in order to achieve a given rate of convergence.

The structure of the present work is as follows: Section 2 is devoted to the presentation of the problem, the clarification of the notion of the solution for a hyperbolic PDE system under boundary sampled-data control, the statement of the main result (Theorem 2.1) and a discussion about the main result. The proof of the main result is provided in Section 3. A simple illustrating example is presented in Section 4. Finally, the concluding remarks are provided in Section 5.

**Notation.** Throughout this paper, we adopt the following notation.

* $\Re_+ := [0,+\infty)$. $Z^+$ denotes the set of all non-negative integers.
* Let $U \subseteq \Re^n$ be a set with non-empty interior and let $\Omega \subseteq \Re$ be a set. By $C^0(U;\Omega)$, we denote the class of continuous mappings on $U$, which take values in $\Omega$. By $C^k(U;\Omega)$, where $k \geq 1$, we denote the class of continuous functions on $U$, which have continuous derivatives of order $k$ on $U$ and take values in $\Omega$.
* Let $I \subseteq \Re$ be an interval. A function $f:I \to \Re$ is called right continuous on $I$ if for every $t \in I$ and $\varepsilon > 0$ there exists $\delta(\varepsilon,t) > 0$ such that for all $\tau \in I$ with $t \leq \tau < t + \delta(\varepsilon,t)$ it holds that $|f(\tau) - f(t)| < \varepsilon$. A function $f:I \to \Re$ is called left continuous on $I$ if for every $t \in I$ and $\varepsilon > 0$ there exists $\delta(\varepsilon,t) > 0$ such that for all $\tau \in I$ with $t \geq \tau > t - \delta(\varepsilon,t)$ it holds that $|f(\tau) - f(t)| < \varepsilon$. A function $f:I \to \Re$ is called piecewise continuous on $I$ if for every compact $K \subseteq I$ the numbers of points $t \in K$ where $f$ is discontinuous is finite and furthermore, for every $t \in I$ all meaningful limits $\lim_{h \to 0^+}(f(t+h))$, $\lim_{h \to 0^+}(f(t-h))$ exist and are finite. Let $a \in \Re$ be a given real number. A function $f:[a,+\infty) \to \Re$ is called right differentiable on $[a,+\infty)$ if for every $t \geq a$ the limit $\lim_{h \to 0^+}\left(h^{-1}(f(t+h) - f(t))\right)$ exists and is finite. A piecewise continuous, right differentiable function

$f:[a,+\infty) \to \Re$ is called piecewise $C^1$ on $[a,+\infty)$ if for every $t > a$ the limit $\lim_{h \to 0^-}\left(h^{-1}(f(t+h) - f(t))\right)$ exists and is finite and furthermore, for every compact $K \subseteq [a,+\infty)$ the numbers of points $t \in K \setminus \{a\}$ where $\lim_{h \to 0^-}\left(h^{-1}(f(t+h) - f(t))\right) \neq \lim_{h \to 0^+}\left(h^{-1}(f(t+h) - f(t))\right)$ is finite.

* Let $x: \Re_+ \times [0,1] \to \Re$ be given. We use the notation $x[t]$ to denote the profile at certain $t \geq 0$, i.e., $(x[t])(z) = x(t,z)$ for all $z \in [0,1]$.

* Let $I \subseteq \Re$ be an interval. $L^2(I)$ denotes the class of measurable functions $f: I \to \Re$ which are square integrable. $L^\infty(I)$ denotes the class of measurable functions $f: I \to \Re$ which are essentially bounded on $I$. $L^\infty_{loc}(I)$ denotes the class of measurable functions $f: I \to \Re$ which are of class $L^\infty(K)$ for every compact $K \subseteq I$.

* We define $p_a(t) := \int_0^t \exp(-a(t-s))ds$ for all $t \geq 0$ and $a \in \Re$. Notice that $p_a(t) := \frac{1 - \exp(-at)}{a}$ for $a \neq 0$ and $p_0(t) := t$.

## 2. Main Result

We consider the control system

$$\frac{\partial y}{\partial t}(t,z) + \frac{\partial y}{\partial z}(t,z) = g(z)y(t,1) + \int_z^1 f(z,s)y(t,s)ds, \text{ for } (t,z) \in \Re_+ \times [0,1] \tag{2.1}$$

$$y(t,0) = u(t) - \int_0^1 p(s)y(t,s)ds, \text{ for } t \geq 0 \tag{2.2}$$

where $g \in C^0([0,1]; \Re)$, $p \in C^1([0,1]; \Re)$, $f \in C^0([0,1]^2; \Re)$ are given functions, $y[t]$ is the state and $u(t)$ is the control input. More specifically, we consider the system under boundary sampled-data control with ZOH:

$$u(t) = u_i, \text{ for } t \in [\tau_i, \tau_{i+1}) \text{ and for all } i \in Z^+ \tag{2.3}$$

where $\{\tau_i \geq 0, i = 0,1,2,...\}$ is an increasing sequence (the sequence of sampling times) with $\tau_0 = 0$, $\lim_{i \to +\infty}(\tau_i) = +\infty$ and $\{u_i \in \Re, i = 0,1,2,...\}$ is the sequence of applied inputs.

Let $y_0 \in C^1([0,1])$ be given. By a solution of the initial value problem (2.1), (2.2), (2.3) with initial condition $y_0$, sampling times $\{\tau_i \geq 0, i = 0,1,2,...\}$ and applied inputs $\{u_i \in \Re, i = 0,1,2,...\}$, we mean a function $y: \Re_+ \times [0,1] \to \Re$, which satisfies the following properties:

a) the mapping $\Re_+ \ni t \to y(t,z) \in \Re$ is piecewise continuous and right differentiable for all $z \in [0,1]$,

b) the mapping $[0,1] \ni z \to y(t,z) \in \Re$ is piecewise continuous and left continuous for all $t \geq 0$,

c) $y: \Re_+ \times [0,1] \to \Re$ is continuous on $B := \{(t,z) \in \Re_+ \times [0,1] : t - z \neq \tau_i, i = 0,1,2,...\}$ and is of class $C^1$ on $A := \{(t,z) \in \Re_+ \times [0,1] : t - z \neq \tau_i, t - z \neq 1 + \tau_i, i = 0,1,2,...\}$,

d) the limit $\lim_{h \to 0^+} \frac{y(t,z) - y(t,z-h)}{h}$ exists and is finite for all $(t,z) \in \Re_+ \times (0,1]$,

e) for every $r > 0$ the mappings $A \cap \{(t,z) \in [0,r] \times [0,1]\} \ni (t,z) \to \frac{\partial y}{\partial t}(t,z) \in \Re$ and $A \cap \{(t,z) \in [0,r] \times [0,1]\} \ni (t,z) \to \frac{\partial y}{\partial z}(t,z) \in \Re$ are bounded,

f) equation (2.1) holds for all $(t,z) \in A$, equations (2.2), (2.3) hold and $y(0,z) = y_0(z)$, for all $z \in (0,1]$.

Theorem 1 in [19] guarantees the existence of functions $k, l \in C^1([0,1]^2)$ such that the Volterra transformation

$$x(t,z) = y(t,z) - \int_z^1 k(z,s) y(t,s) ds, \text{ for all } (t,z) \in \Re_+ \times [0,1] \quad (2.4)$$

with inverse

$$y(t,z) = x(t,z) + \int_z^1 l(z,s) x(t,s) ds, \text{ for all } (t,z) \in \Re_+ \times [0,1] \quad (2.5)$$

maps the solutions of (2.1), (2.2) to the solutions of

$$\frac{\partial x}{\partial t}(t,z) + \frac{\partial x}{\partial z}(t,z) = 0, \text{ for } (t,z) \in \Re_+ \times [0,1] \quad (2.6)$$

$$x(t,0) = u(t) - \int_0^1 (p(s) + k(0,s)) y(t,s) ds, \text{ for } t \geq 0, \quad (2.7)$$

More specifically, the kernels $k, l \in C^1([0,1]^2)$ satisfy the following equations for all $z, s \in [0,1]$:

$$\begin{aligned} l(z,1) &= -g(z) \\ \frac{\partial l}{\partial z}(z,s) + \frac{\partial l}{\partial s}(z,s) &= f(z,s) + \int_z^s f(z,w) l(w,s) dw \end{aligned} \quad (2.8)$$

$$\begin{aligned} k(z,1) &= -g(z) + \int_z^1 k(z,s) g(s) ds \\ \frac{\partial k}{\partial z}(z,s) + \frac{\partial k}{\partial s}(z,s) &= f(z,s) - \int_z^s k(z,w) f(w,s) dw \end{aligned} \quad (2.9)$$

$$l(z,s) = k(z,s) + \int_z^s l(z,w) k(w,s) dw \quad (2.10)$$

Using the emulation sampled-data design based on the continuous time feedback law $u(t) = \int_0^1 (p(s) + k(0,s)) y(t,s) ds$ for the control system (2.1), (2.2), (2.3), we are in a position to establish the following stabilization result.

**Theorem 2.1:** *For every $\sigma > 0$ there exist constants $T, G > 0$ with the following property: for every $y_0 \in C^1([0,1])$ and for every increasing sequence $\{\tau_i \geq 0, i = 0,1,2,...\}$ with $\tau_0 = 0$, $\sup_{i \geq 0}(\tau_{i+1} - \tau_i) \leq T$ and $\lim_{i \to +\infty}(\tau_i) = +\infty$ the initial value problem of the closed-loop system (2.1), (2.2), (2.3) with*

$$u_i = \int_0^1 (p(s) + k(0,s)) y(\tau_i, s) ds, \text{ for all } i \in Z^+ \quad (2.11)$$

*and initial condition $y_0$ has a unique solution, which also satisfies the following estimate:*

$$\|y[t]\|_\infty \leq G \exp(-\sigma t) \|y_0\|_\infty, \text{ for all } t \geq 0 \quad (2.12)$$

The proof of Theorem 2.1 is provided in the following section. The ideas behind the proof of Theorem 2.1 are described next.

**Idea 1:** In order to establish estimate (2.12) for $y[t]$, it suffices to establish a similar estimate for $x[t]$, where $x$ is the solution of (2.6), (2.7), (2.3), (2.12) with initial condition $x_0 \in C^1([0,1])$ that satisfies $x_0(z) = y_0(z) - \int_z^1 k(z,s) y_0(s) ds$ for all $z \in [0,1]$ (see (2.4)).

**Idea 2:** The solution of (2.6), (2.7), (2.3), (2.12) can be determined without knowledge of $y[t]$. Indeed, using (2.5), (2.7) and (2.11), we get:

$$x(t,0) = u(t) - \int_0^1 (p(s) + k(0,s)) x(t,s) ds - \int_0^1 \int_s^1 (p(s) + k(0,s)) l(s,w) x(t,w) dw ds = u(t) - \int_0^1 \tilde{k}(s) x(t,s) ds, \text{ for } t \geq 0,$$
(2.13)

$$u_i = \int_0^1 \tilde{k}(s) x(\tau_i, s) ds, \text{ for } t \geq 0,$$
(2.14)

where

$$\tilde{k}(s) := p(s) + k(0,s) + \int_0^s (p(w) + k(0,w)) l(w,s) dw, \text{ for all } s \in [0,1]$$
(2.15)

So, we need to establish an exponential stability estimate (similar to (2.12)) for the solution $x[t]$ of (2.3), (2.6), (2.13), (2.14).

**Idea 3:** Since the solution $x[t]$ of (2.3), (2.6), (2.13), (2.14) will be (in general) discontinuous, we need to work with the mild solution of (2.3), (2.6), (2.13), (2.14) (see [16] for the notion of the mild solution). Moreover, the mild solution of (2.3), (2.6), (2.13), (2.14) will be constructed by an IDE (see [16]). Consequently, we need to construct solutions for a particular IDE and we need to show exponential stability estimates for the solutions of this particular IDE. To this end, we are using the fact that for every $\sigma > 0$, $x_0 \in L^\infty((0,1])$, $v \in L^\infty_{loc}(\Re_+)$, the unique mild solution of the evolution problem

$$\frac{\partial x}{\partial t}(t,z) + \frac{\partial x}{\partial z}(t,z) = 0, \text{ for } (t,z) \in \Re_+ \times [0,1]$$
(2.16)

$$x(t,0) = v(t), \text{ for } t \geq 0$$
(2.17)

with initial condition $x[0] = x_0$ satisfies the estimate:

$$\|x[t]\|_\infty \leq \exp(-\sigma(t-1)) \|x_0\|_\infty + \exp(\sigma) \sup_{\max(0,t-1) \leq s \leq t} \left( |v(s)| \exp(-\sigma(t-s)) \right)$$
(2.18)

The above fact is a direct consequence of the notion of the mild solution given in [16] and the fact that the mild solution of the evolution problem (2.16), (2.17) is given by the formula:

$$x(t,z) = \begin{cases} v(t-z) & \text{for } t \geq z \\ x_0(z-t) & \text{for } t < z \end{cases}$$

**Remark 2.2:** A comparison between the sampled-data stabilization of parabolic PDEs given in [18] and Theorem 2.1 is fruitful.

**(i)** For parabolic PDEs both main results in [18] (Theorem 3.1 and Theorem 3.2) provide exponential stability estimates in the $L^2$ norm, while Theorem 2.1 provides an exponential stability estimate (estimate (2.12)) in the $L^\infty$ norm.

**(ii)** The proofs of both main results in [18] are very different from the proof of Theorem 2.1. While both proofs of the main results in [18] exploit an eigenfunction expansion procedure, the proof of Theorem 2.1 relies on the representation of the solution of the closed-loop system (2.1), (2.2), (2.3) with (2.11) by means of an IDE.

**(iii)** Theorem 2.1 as well as both main results in [18] guarantee robustness with respect to perturbations of the sampling schedule: the exponential stability estimate (2.12) holds for every increasing sequence $\{\tau_i \geq 0, i = 0,1,2,...\}$ with $\tau_0 = 0$, $\sup_{i \geq 0}(\tau_{i+1} - \tau_i) \leq T$, $\lim_{i \to +\infty}(\tau_i) = +\infty$.

**(iv)** Finally, there is an additional important difference between the parabolic and the hyperbolic case. In the hyperbolic case (Theorem 2.1), by selecting a sufficiently small maximum allowable sampling period $T > 0$ we can achieve an arbitrary rate of convergence $\sigma > 0$. This is not possible in the parabolic case. This important difference can be explained by the fact that the nominal continuous-time feedback law $u(t) = \int_0^1 (p(s) + k(0,s))y(t,s)ds$ achieves finite-time stability in the hyperbolic case, while the nominal continuous-time feedback law in the parabolic case achieves exponential stability. The proof of Theorem 2.1 provides an estimate of how small the maximum allowable sampling period $T > 0$ must be in order to achieve a given rate of convergence $\sigma > 0$: the inequality

$$Mp_a(T)\exp(\sigma(1+T)) < 1 \tag{2.19}$$

must be satisfied, where

$$M := \left|\tilde{k}(1)\right| + \int_0^1 \left|\frac{d\tilde{k}}{ds}(s)\right|ds, \quad a := \tilde{k}(0) \tag{2.20}$$

Formulas (2.19), (2.20) can be used by the control practitioner in order to obtain a qualitative estimate for the maximum allowable sampling period $T > 0$. However, it should be noticed that inequality (2.19) is conservative and larger values for the maximum allowable sampling period $T > 0$ can be used in practice. Figure 1 shows the dependence of the rate of convergence $\sigma > 0$ on the maximum allowable sampling period $T > 0$.

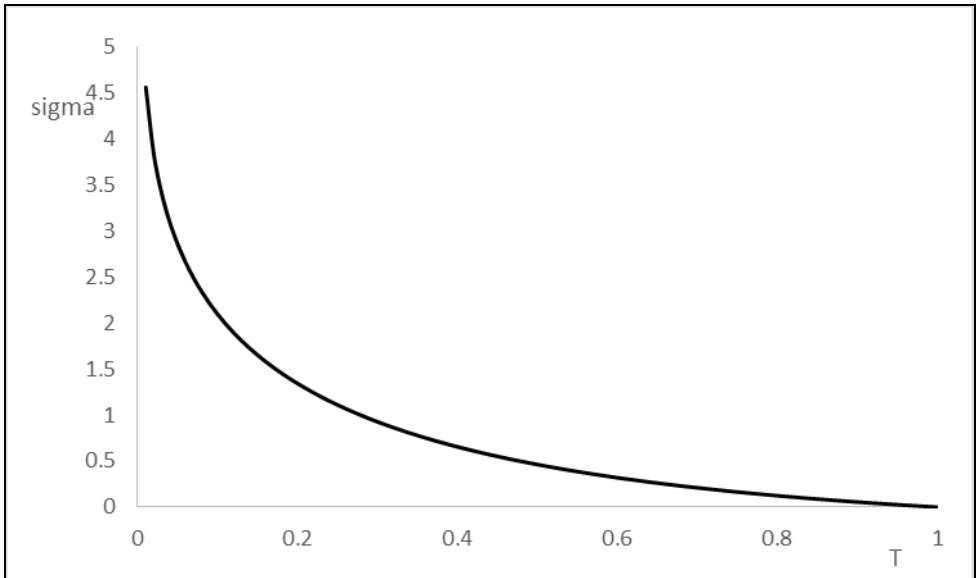

**Fig. 1:** The dependence of the rate of convergence $\sigma > 0$ on the maximum allowable sampling period $T > 0$, as predicted by inequality (2.19) for $a = 0$, $M = 1$.

# 3. Proof of Theorem 2.1

The rigorous presentation of the proof of Theorem 2.1 follows the aforementioned ideas and relies on certain facts. The first fact is a technical fact that allows us to construct the solution of a specific IDE on an infinite interval.

**Fact I:** *Let $\tilde{k} \in C^1([0,1])$ be a given function. Then for every $u_i \in \Re$, $\tau_i \geq 0$ and for every $v_0 \in L^\infty([\tau_i -1, \tau_i))$ being right continuous and piecewise continuous, the initial value problem of the IDE*

$$v(t) = u_i - \int_0^1 \tilde{k}(s) v(t-s) ds \tag{3.1}$$

*for $t \geq \tau_i$, with initial condition $v(s) = v_0(s)$ for $s \in [\tau_i -1, \tau_i)$, has a unique right-continuous solution $v:[\tau_i -1, +\infty) \to \Re$, which is continuous on $(\tau_i, +\infty)$. Moreover, $v(t)$ is piecewise $C^1$ and right differentiable on $[\tau_i, +\infty)$. Finally, if $B = \{t_1, t_2, ..., t_m\} \subseteq [\tau_i -1, \tau_i]$ is the (finite) set of points of discontinuities of $v(t)$, then $v(t)$ is $C^1$ on $(\tau_i, +\infty) \setminus B'$ and satisfies for all $t \in (\tau_i, +\infty) \setminus B'$ the differential equation*

$$\dot{v}(t) = -\tilde{k}(0) v(t) + \tilde{k}(1) v(t-1) - \int_0^1 \frac{d\tilde{k}}{ds}(s) v(t-s) ds \tag{3.2}$$

*where $B' = \{1+t_1, 1+t_2, ..., 1+t_m\}$.*

**Proof of Fact I:** Let $u_i \in \Re$, $\tau_i \geq 0$ and $v_0 \in L^\infty([\tau_i -1, \tau_i))$ be given (arbitrary). Theorem 2.1 in [16] guarantees that there exists $t_{max} \in (0, +\infty]$ and a unique mapping $v \in L^\infty_{loc}([\tau_i -1, \tau_i + t_{max}))$ such that $v(s) = v_0(s)$ for $s \in [\tau_i -1, \tau_i)$ and (3.1) holds for $t \in [\tau_i, \tau_i + t_{max})$ a.e.. It follows from (3.1) and its equivalent form

$$v(t) = u_i - \int_{t-1}^t \tilde{k}(t-s) v(s) ds, \tag{3.3}$$

that $v \in L^\infty_{loc}([\tau_i -1, \tau_i + t_{max}))$ is a right continuous and piecewise continuous function satisfying (3.1) for all $t \in [\tau_i, \tau_i + t_{max})$. Moreover, $v(t)$ is continuous on $(\tau_i, \tau_i + t_{max})$.

We will show next that $t_{max} = +\infty$. Indeed, if $t_{max} < +\infty$, then it will follow from Theorem 2.1 in [16] that $\lim_{t \to t_{max}^-} \sup_{\tau_i \leq s < t} (|v(s)|) = +\infty$. Let $R, G > 0$ be positive constants so that

$$R \geq \int_0^1 |\tilde{k}(s)| ds, \quad G \geq \max_{s \in [0,1]} (|\tilde{k}(s)|) \tag{3.4}$$

Equation (3.3) implies that for all $t \in [\tau_i, \tau_i + \min(1, t_{max}/2)]$ the following equation holds:

$$v(t) = u_i - \int_{\tau_i}^t \tilde{k}(t-s) v(s) ds - \int_{t-1}^{\tau_i} \tilde{k}(t-s) v_0(s) ds \tag{3.5}$$

Using inequalities (3.4), (3.4), we obtain from (3.5) for all $t \in [\tau_i, \tau_i + \min(1, t_{max}/2)]$:

$$|v(t)| \leq |u_i| + G(t-\tau_i) \sup_{\tau_i \leq s \leq t}(|v(s)|) + R \sup_{\tau_i - 1 \leq s < \tau_i}(|v_0(s)|) \tag{3.6}$$

Define:

$$h := \frac{\min(1, t_{\max})}{2(1+G)} \tag{3.7}$$

It follows from (3.6) and (3.7) that the following estimate holds:

$$\sup_{\tau_i \leq s \leq t}(|v(s)|) \leq 2\left(|u_i| + R \sup_{\tau_i - 1 \leq s < \tau_i}(|v_0(s)|)\right) \text{ for } t \in [\tau_i, \tau_i + h] \tag{3.8}$$

Without loss of generality, we may assume that $R \geq 1$. Then we obtain from (3.8) and the fact that $v(s) = v_0(s)$ for $s \in [\tau_i - 1, \tau_i)$ the following estimate:

$$\sup_{\tau_i - 1 \leq s < \tau_i + h}(|v(s)|) \leq 2\left(|u_i| + R \sup_{\tau_i - 1 \leq s < \tau_i}(|v_0(s)|)\right) \tag{3.9}$$

Applying (3.9) repeatedly, we obtain the following estimate for all integers $m \geq 1$:

$$\sup_{\tau_i - 1 \leq s < \tau_i + T_m}(|v(s)|) \leq (2R)^m \sup_{\tau_i - 1 \leq s < \tau_i}(|v_0(s)|) + 2|u_i|\frac{(2R)^{m-1} - 1}{2R - 1} \tag{3.10}$$

where $T_m = \min(mh, t_{\max})$. Estimate (3.10) contradicts the fact that $\lim_{t \to t_{\max}^-} \sup_{\tau_i \leq s < t}(|v(s)|) = +\infty$. Therefore, we must have $t_{\max} = +\infty$.

It follows from (3.3) that $v$ is right differentiable on $[\tau_i, +\infty)$ with

$$D^+ v(t) = -\tilde{k}(0)v(t) + \tilde{k}(1)v(t-1) - \int_0^1 \frac{d\tilde{k}}{ds}(s) v(t-s) ds \tag{3.11}$$

for all $t \geq \tau_i$, where $D^+ v(t)$ denotes the right derivative of $v$ at $t$. It follows from (3.11) that $D^+ v(t)$ is continuous on $(\tau_i, +\infty) \setminus B'$ with $B' = \{1+t_1, 1+t_2, \ldots, 1+t_m\}$. Corollary 1.2 on page 43 in [25] implies that $v(t)$ is $C^1$ on $(\tau_i, +\infty) \setminus B'$ and satisfies (3.2) on $(\tau_i, +\infty) \setminus B'$. Finally, (3.2) implies that $v(t)$ is piecewise $C^1$ on $[\tau_i, +\infty)$ (since $\dot{v}(t)$ has finite left and right limits everywhere). The proof of the fact is complete. ◁

Applying Fact I repeatedly, we obtain the following fact.

**Fact II:** *Let $\tilde{k} \in C^1([0,1])$ be a given function. Then for every increasing sequence $\{\tau_i \geq 0, i = 0,1,2,\ldots\}$ with $\tau_0 = 0$, $\lim_{i \to +\infty}(\tau_i) = +\infty$, for every sequence $\{u_i \in \Re, i = 0,1,2,\ldots\}$ and for every $v_0 \in L^\infty([-1,0))$ being right continuous and piecewise continuous, the IDE (3.1) for $t \in [\tau_i, \tau_{i+1})$ and $i \in Z^+$ with initial condition $v(s) = v_0(s)$ for $s \in [-1,0)$, has a unique right-continuous solution $v:[-1,+\infty) \to \Re$, which is continuous on $\Re_+ \setminus \{\tau_i \geq 0, i = 0,1,2,\ldots\}$. Moreover, $v(t)$ is piecewise $C^1$ and right differentiable on $\Re_+$. Finally, if $B = \{t_1, t_2, \ldots, t_m\} \subseteq [-1,0)$ is the (finite) set of points of discontinuities of $v_0$, then $v(t)$ is $C^1$ on $\Re_+ \setminus (B' \cup \{1+\tau_i \geq 0: i = 0,1,2,\ldots\} \cup \{\tau_i \geq 0: i = 0,1,2,\ldots\})$ and satisfies for all $t \in \Re_+ \setminus (B' \cup \{1+\tau_i \geq 0: i = 0,1,2,\ldots\} \cup \{1+\tau_i \geq 0: i = 0,1,2,\ldots\})$ the differential equation (3.2), where $B' = \{1+t_1, 1+t_2, \ldots, 1+t_m\}$.*

The following fact provides an exponential stability estimate for the IDE (3.1) with sampled-data control. The stability analysis follows a small-gain technique.

**Fact III:** *Let $\tilde{k} \in C^1([0,1])$ be a given function. Let $\sigma, T > 0$ be constants satisfying inequality (2.19), where the constants $a \in \Re$, $M \geq 0$ are given by (2.20). Then for every increasing sequence $\{\tau_i \geq 0, i = 0,1,2,...\}$ with $\tau_0 = 0$, $\sup_{i \geq 0}(\tau_{i+1} - \tau_i) \leq T$, $\lim_{i \to +\infty}(\tau_i) = +\infty$ and for every $v_0 \in L^\infty([-1,0))$ being right continuous, the unique solution $v(t)$ of the IDE*

$$v(t) = u(t) - \int_0^1 \tilde{k}(s) v(t-s) ds, \text{ for } t \geq 0 \tag{3.12}$$

*with*

$$u(t) = u_i = \int_0^1 \tilde{k}(s) v(\tau_i - s) ds, \text{ for all } t \in [\tau_i, \tau_{i+1}) \text{ and for all } i \in Z^+, \tag{3.13}$$

*and initial condition $v(s) = v_0(s)$ for $s \in [-1,0)$, satisfies the estimate*

$$|v(t)| \exp(\sigma t) \leq \sup_{-1 \leq s < 0}(|v_0(s)|), \text{ for all } t \geq -1 \tag{3.14}$$

**Proof of Fact III:** Fact II guarantees that for every integer $i \geq 0$, differential equation (3.2) holds for $t \in (\tau_i, \tau_{i+1})$ a.e.. It follows from (3.12) and (3.13) that $v(\tau_i) = 0$ for every integer $i \geq 0$. Using this fact, right continuity of $v:[-1,+\infty) \to \Re$ and absolute continuity of $v(t)$ on $(\tau_i, \tau_{i+1})$, we get the following formula that holds for all $t \in [\tau_i, \tau_{i+1})$:

$$v(t) = \exp(-a(t-\tau_i))v(\tau_i) + \int_{\tau_i}^t \exp(-a(t-s))\left(\tilde{k}(1)v(s-1) - \int_0^1 \frac{d\tilde{k}}{ds}(l)v(s-l)dl\right)ds \tag{3.15}$$

where $a = \tilde{k}(0)$. Using (3.15) and definitions (2.20) we obtain for all $t \in [\tau_i, \tau_{i+1})$ and every integer $i \geq 0$:

$$|v(t)| \leq p_a(t-\tau_i) M \sup_{\tau_i - 1 \leq s \leq t}(|v(s)|) \tag{3.16}$$

A direct consequence of (3.16) is the following estimate that holds for all $t \in [\tau_i, \tau_{i+1})$:

$$|v(t)| \exp(\sigma t) \leq M \, p_a(t-\tau_i) \exp(\sigma(t-\tau_i) + \sigma) \sup_{\tau_i - 1 \leq s \leq t}(|v(s)| \exp(\sigma s)) \tag{3.17}$$

It follows from (3.17), the fact that $\sup_{i \geq 0}(\tau_{i+1} - \tau_i) \leq T$ and the fact that the mapping $\Re_+ \ni t \to p_a(t)$ is non-decreasing for all $a \in \Re$, that the following estimate holds for all $t \in [\tau_i, \tau_{i+1})$:

$$|v(t)| \exp(\sigma t) \leq M \, p_a(T) \exp(\sigma(1+T)) \sup_{-1 \leq s \leq t}(|v(s)| \exp(\sigma s)) \tag{3.18}$$

At this point, it should be noticed that estimate (3.18) holds for all $t \geq 0$, since there is no dependence on the integer $i \geq 0$. Consequently, we obtain from (3.18) the following estimate that holds for all $t \geq 0$:

$$\sup_{0 \leq s \leq t}(|v(s)| \exp(\sigma s)) \leq M p_a(T) \exp(\sigma(1+T)) \max\left(\sup_{-1 \leq s < 0}(|v(s)| \exp(\sigma s)), \sup_{0 \leq s \leq t}(|v(s)| \exp(\sigma s))\right) \tag{3.19}$$

Distinguishing the cases $\max\left(\sup_{-1\leq s<0}(|v(s)|\exp(\sigma s)), \sup_{0\leq s\leq t}(|v(s)|\exp(\sigma s))\right) = \sup_{0\leq s\leq t}(|v(s)|\exp(\sigma s))$ and $\max\left(\sup_{-1\leq s<0}(|v(s)|\exp(\sigma s)), \sup_{0\leq s\leq t}(|v(s)|\exp(\sigma s))\right) = \sup_{-1\leq s<0}(|v(s)|\exp(\sigma s))$ in conjunction with inequality (2.19), we obtain from (3.19) for all $t \geq 0$:

$$\sup_{0\leq s\leq t}(|v(s)|\exp(\sigma s)) \leq TM \exp(\sigma(1+T)) \sup_{-1\leq s<0}(|v(s)|\exp(\sigma s)) \quad (3.20)$$

Inequality (3.20) in conjunction with the fact that $v(s) = v_0(s)$ for $s \in [-1,0)$ implies estimate (3.14). The proof of the fact is complete. ◁

We are now ready to provide the proof of Theorem 2.1.

**Proof of Theorem 2.1:** Define $\tilde{k} \in C^1([0,1])$ by means of (2.15). Let $\sigma > 0$ be a given constant and let $T > 0$ be a sufficiently small constant so that (2.19) holds.

Let (arbitrary) $y_0 \in C^1([0,1])$ be given and let an (arbitrary) increasing sequence $\{\tau_i \geq 0, i = 0,1,2,...\}$ with $\tau_0 = 0$, $\sup_{i\geq 0}(\tau_{i+1} - \tau_i) \leq T$ and $\lim_{i\to+\infty}(\tau_i) = +\infty$ be also given. Define

$$v_0(-z) := y_0(z) - \int_z^1 k(z,s) y_0(s) ds \text{ for all } z \in (0,1] \quad (3.21)$$

It follows from Fact II and definition (3.21) (which guarantees that $v_0 \in C^1([-1,0))$ with finite limits for $v_0$ and $\dot{v}_0$ as $t \to 0^-$) that the IDE (3.12), (3.13) with initial condition $v(s) = v_0(s)$ for $s \in [-1,0)$, has a unique right-continuous solution $v:[-1,+\infty) \to \Re$, which is continuous on $\Re_+ \setminus \{\tau_i \geq 0, i = 0,1,2,...\}$. Moreover, $v(t)$ is piecewise $C^1$ and right differentiable on $\Re_+$. Finally, $v(t)$ is $C^1$ on $\Re_+ \setminus (\{1+\tau_i \geq 0: i = 0,1,2,...\} \cup \{\tau_i \geq 0: i = 0,1,2,...\})$ and satisfies for all $t \in \Re_+ \setminus (\{1+\tau_i \geq 0: i = 0,1,2,...\} \cup \{\tau_i \geq 0: i = 0,1,2,...\})$ the differential equation (3.2).

It follows from Fact III that estimate (3.14) holds.

Next, define

$$y(t,z) := v(t-z) + \int_z^1 l(z,s) v(t-s) ds = v(t-z) + \int_{t-1}^{t-z} l(z, t-w) v(w) dw, \text{ for } (t,z) \in \Re_+ \times [0,1] \quad (3.22)$$

It follows from (2.10), (3.21) and (3.22) that $y(0,z) = y_0(z)$, for all $z \in (0,1]$. Moreover, the regularity properties of guarantee that the function $y: \Re_+ \times [0,1] \to \Re$ satisfies the following properties:
  a) the mapping $\Re_+ \ni t \to y(t,z) \in \Re$ is piecewise continuous and right differentiable for all $z \in [0,1]$,
  b) the mapping $[0,1] \ni z \to y(t,z) \in \Re$ is piecewise continuous and left continuous for all $t \geq 0$,
  c) $y: \Re_+ \times [0,1] \to \Re$ is continuous on $B := \{(t,z) \in \Re_+ \times [0,1] : t - z \neq \tau_i, i = 0,1,2,...\}$ and is of class $C^1$ on $A := \{(t,z) \in \Re_+ \times [0,1] : t - z \neq \tau_i, t - z \neq 1 + \tau_i, i = 0,1,2,...\}$,
  d) the limit $\lim_{h \to 0^+} \frac{y(t,z) - y(t,z-h)}{h}$ exists and is finite for all $(t,z) \in \Re_+ \times (0,1]$,

e) for every $r > 0$ the mappings $A \cap \{(t,z) \in [0,r] \times [0,1]\} \ni (t,z) \to \frac{\partial y}{\partial t}(t,z) \in \Re$ and $A \cap \{(t,z) \in [0,r] \times [0,1]\} \ni (t,z) \to \frac{\partial y}{\partial z}(t,z) \in \Re$ are bounded.

Using (2.8) and definition (3.22) that equation (2.1) holds for all $(t,z) \in A$. Finally, it follows from (3.12), (3.13) and definitions (2.15), (3.22) that equations (2.2), (2.3), (2.11) hold.

Definition (3.21) implies the existence of a constant $P_1 > 0$ such that

$$\sup_{-1 \le s < 0} (|v_0(s)|) \le P_1 \|y_0\|_\infty \tag{3.23}$$

Notice that the constant $P_1 > 0$ is independent of $y_0 \in C^1([0,1])$. Moreover, definition (3.22) implies the existence of a constant $P_2 > 0$ such that

$$\|y[t]\|_\infty \le P_2 \sup_{-1 \le s \le 0} (|v(t+s)|), \text{ for all } t \ge 0 \tag{3.24}$$

Again, notice that the constant $P_2 > 0$ is independent of $v : [-1, +\infty) \to \Re$. Estimate (2.12) is a direct consequence of estimates (3.14), (3.23) and (3.24). The proof is complete. ◁

## 4. Illustrative Example

We consider the following control system

$$\frac{\partial y}{\partial t}(t,z) + \frac{\partial y}{\partial z}(t,z) = A\exp(rz) y(t,1), \text{ for } (t,z) \in \Re_+ \times [0,1] \tag{4.1}$$

$$y(t,0) = u(t), \text{ for } t \ge 0 \tag{4.2}$$

where $A, r > 0$ are positive constants, $y[t]$ is the state and $u(t)$ is the control input. System (4.1), (4.2) is a system of the form (2.1), (2.2) with $g(z) = A\exp(rz)$, $f(z,s) \equiv 0$ and $p(s) \equiv 0$. A stability analysis of the open-loop system (4.1), (4.2) with $u(t) \equiv 0$ may be performed by converting the system (4.1), (4.2) to a system of IDEs and applying Theorem 3.7 in [16]. Indeed, system (4.1), (4.2) with $u(t) \equiv 0$ is equivalent to the following system of IDEs (see [16])

$$v(t) = 0$$

$$p(t) = v(t-1) + A \int_0^1 \exp(rw) p(t-1+w) dw$$

and applying Theorem 3.7 in [16] with $N = 1$, $c = 1$, $W(p,v) := K|v| + |p|$ and arbitrary $K > 1$, we can conclude that the open-loop system (4.1), (4.2) with $u(t) \equiv 0$ is Globally Asymptotically Stable provided that

$$A < \frac{r}{\exp(r) - 1} \tag{4.3}$$

On the other hand, when condition (4.3) is violated, then there are solutions which do not converge to zero: the functions $x(t,z) = \exp(\lambda(t+1-z)) \frac{\exp((r+\lambda)z) - 1}{\exp(r+\lambda) - 1}$, where $\lambda \ge 0$ is the unique solution of the equation $A = \frac{(r+\lambda)\exp(\lambda)}{\exp(r+\lambda)-1}$ (guaranteed to exist when (4.3) does not hold) are solutions of the open-loop system (4.1), (4.2) with $u(t) \equiv 0$. Therefore, the stability condition (4.3) is sharp.

When (4.3) does not hold, there is need for control. The case of applying a continuous-time feedback law was considered in Example 2.1 in [19] (with different space orientation), where finite-time stabilization was achieved. Here, we consider the system under boundary sampled-data control with ZOH:

$$u(t) = u_i, \text{ for } t \in [\tau_i, \tau_{i+1}) \text{ and for all } i \in Z^+ \tag{4.4}$$

where $\{\tau_i \geq 0, i = 0,1,2,...\}$ is an increasing sequence (the sequence of sampling times) with $\tau_0 = 0$, $\lim_{i \to +\infty}(\tau_i) = +\infty$ and $\{u_i \in \Re, i = 0,1,2,...\}$ is the sequence of applied inputs.

For this problem, the kernels $k, l \in C^1([0,1]^2)$ that satisfy (2.8), (2.9) and (2.10) are given by the following equations for all $z, s \in [0,1]$:

$$k(z,s) := -A\exp(r)\exp\big((r - A\exp(r))(z-s)\big) \tag{4.5}$$

$$l(z,s) := -A\exp(r(z-s+1)) \tag{4.6}$$

It follows from Theorem 2.1 that for every $\sigma > 0$ there exist constants $T, G > 0$ with the following property: for every $y_0 \in C^1([0,1])$ and for every increasing sequence $\{\tau_i \geq 0, i = 0,1,2,...\}$ with $\tau_0 = 0$, $\sup_{i \geq 0}(\tau_{i+1} - \tau_i) \leq T$ and $\lim_{i \to +\infty}(\tau_i) = +\infty$ the initial value problem of the closed-loop system (4.1), (4.2), (4.4) with

$$u_i = -A\exp(r)\int_0^1 \exp\big((A\exp(r) - r)s\big)y(\tau_i, s)ds, \text{ for all } i \in Z^+ \tag{4.7}$$

and initial condition $y_0$ has a unique solution, which also satisfies estimate (2.12). The proof of Theorem 2.1 guarantees that for every $\sigma > 0$, the maximum allowable sampling period $T > 0$ for which estimate (2.12) holds may be selected so as to satisfy condition (2.19), where $M \geq 0$ and $a \in \Re$ are given by (2.20). For this example, we have:

$$\tilde{k}(s) := -A\exp(r)\exp(-rs), \text{ for all } s \in [0,1] \tag{4.8}$$

and

$$M := A\exp(r), \quad a := -A\exp(r) \tag{4.9}$$

Formulas (4.8), (4.9) can be used in order to obtain a qualitative estimate for the maximum allowable sampling period $T > 0$ from inequality (2.19). However, since inequality (2.19) is conservative, larger values for the maximum allowable sampling period $T > 0$ can be used in practice.

Figure 2 shows the evolution of $y(t,1)$ for the closed-loop system (4.1), (4.2), (4.4) with (4.7), $\tau_i = iT$ for $i \in Z^+$, $T = 0.1$, $r = A = 1$ and initial condition $y_0(z) = -\frac{1}{2}(\exp(z) - e + 1)$ for $z \in [0,1]$. Figure 3 shows the evolution of the input $u(t)$ for the same case. Notice that the open-loop system (4.1), (4.2) with $u(t) \equiv 0$ is unstable, since (4.3) does not hold. For this particular case, the convergence rate predicted from condition (2.19), (4.8), (4.9) satisfies $\sigma < 1.058$.

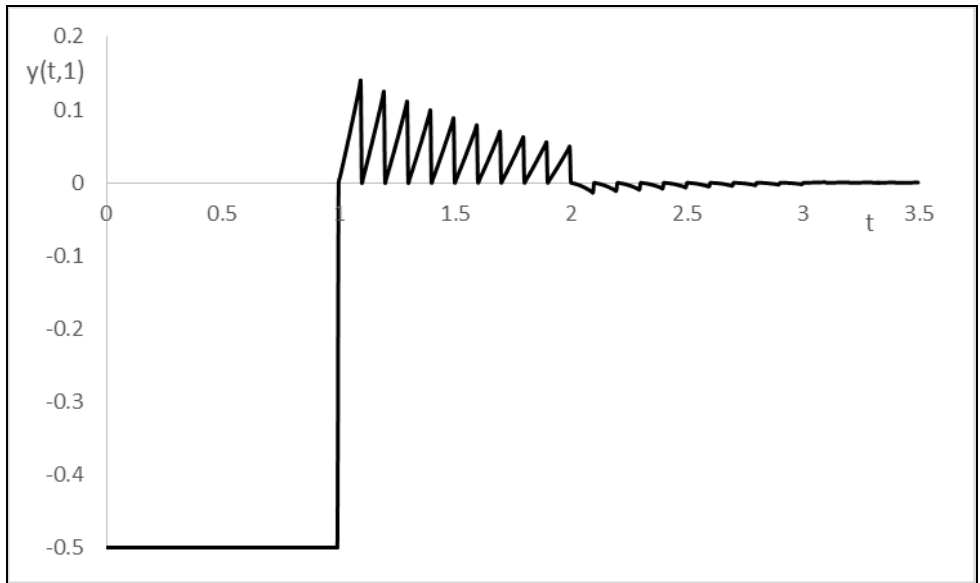

**Fig. 2:** The evolution of $y(t,1)$ for the closed-loop system (4.1), (4.2), (4.4) with (4.7), $\tau_i = iT$ for $i \in Z^+$, $T = 0.1$, $r = A = 1$ and initial condition $y_0(z) = -\frac{1}{2}(\exp(z) - e + 1)$ for $z \in [0,1]$.

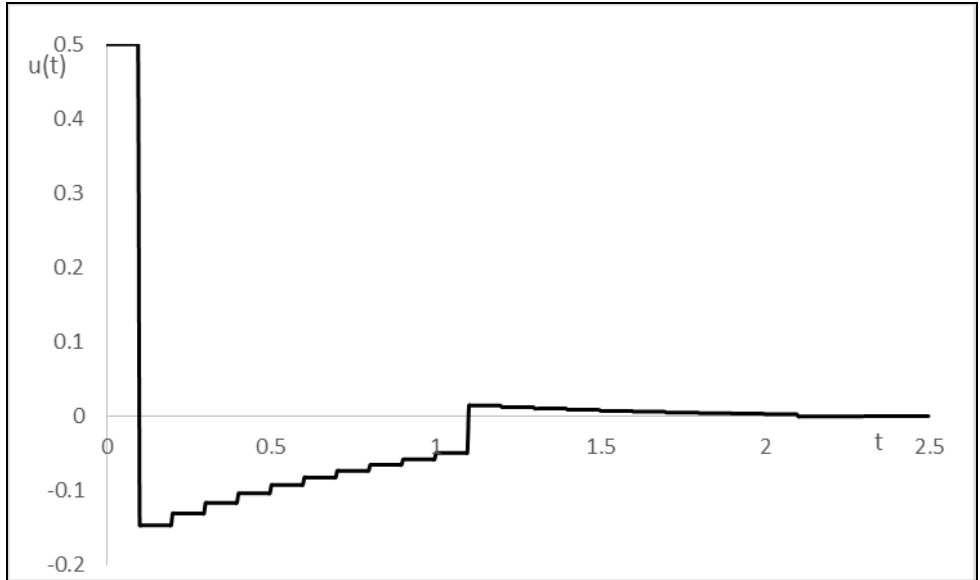

**Fig. 3:** The evolution of the input $u(t)$ for the closed-loop system (4.1), (4.2), (4.4) with (4.7), $\tau_i = iT$ for $i \in Z^+$, $T = 0.1$, $r = A = 1$ and initial condition $y_0(z) = -\frac{1}{2}(\exp(z) - e + 1)$ for $z \in [0,1]$.

## 5. Concluding Remarks

The paper provides results for the application of boundary feedback control with ZOH to 1-D linear, first-order, hyperbolic systems with non-local terms on bounded domains. The control design is based on the continuous-time, boundary feedback, designed by means of backstepping in [19], which guarantees finite-time stability for the corresponding closed-loop system. The application of the aforementioned feedback design with ZOH guarantees closed-loop exponential stability, provided that the sampling period is sufficiently small. It is also shown that, contrary to the parabolic case, a smaller sampling period implies a faster convergence rate with no upper bound for the achieved convergence rate. The obtained results provide stability estimates for the sup-norm of the state and robustness with respect to perturbations of the sampling schedule is guaranteed.

Future work may involve the development of boundary feedback designs that are capable to handle simultaneous sampling in space and time. To this purpose, sampled-data observers for hyperbolic PDEs with non-local terms must be developed.